\documentclass[12pt]{article}

\usepackage{amssymb}
\usepackage{array} % for better arrays (eg matrices) in maths
\usepackage{amsthm}

\theoremstyle{theorem}
\newtheorem*{lemma}{Lemma}

\theoremstyle{definition}
\newtheorem*{definition}{Definition}

\begin{document}
\title{Groups With At Most Twelve Subgroups}
\author{Michael C. Slattery\\
Marquette University, Milwaukee WI 53201-1881\\
\texttt{michael.slattery@marquette.edu}}

\maketitle

In \cite{Slat} I showed that the finite groups with a specified number of subgroups
can always be described as a finite list of similarity classes.  Neil Sloane suggested
that I submit the corresponding sequence, the number of similarity classes with
$n$ subgroups, to his Online Encyclopedia of Integer Sequences \cite{OEIS}.
I thought this would be a quick calculation until I discovered that
my ``example" case in the paper ($n = 6$) was wrong (as noted in \cite{Slat2}).
I realized that producing a reliable count required a full proof.  This note contains
the proof behind my computation of the first 12 terms of the sequence.

For completeness, we include a section of \cite{Slat}:

\begin{quotation}
Now, we saw [above] that cyclic Sylow subgroups of $G$ which
are direct factors allow many non-isomorphic groups to have the
same number of subgroups.  [This lemma] tells us that these are
precisely the cyclic Sylow subgroups which lie in the center of
$G$.  Consequently, if $p_1, \ldots, p_c$ are the primes which
divide $|G|$ such that a Sylow $p_i$-subgroup is cyclic and
central, then we can write
$$
G = P_1 \times P_2 \times \cdots \times P_c \times O_{\pi'}(G)
$$
where each $P_i$ is a Sylow $p_i$-subgroup, the set
$\pi = \{p_1, \ldots, p_c\}$ and $O_{\pi'}(G)$ is the largest normal subgroup
of $G$ with order not divisible by any prime in $\pi$.  We will write
$\widetilde{G} = O_{\pi'}(G)$.  In other words, $\widetilde{G}$
is
the unique subgroup of $G$ left after factoring out the cyclic,
central Sylow subgroups.  It is the part of $G$ that we can hope
to restrict in terms of the number of subgroups.  On the other
hand, if we substituted a different prime (relatively prime to
$|G|$) for any one of the $p_i$, it would not affect the number
of subgroups of $G$.  Thus, we are led to define the following
equivalence relation (which we will call ``similar'' for this
note).

\begin{definition}
Let $G$ and $H$ be two finite groups.
Write $G$ as
a
product of cyclic central Sylow subgroups and $\widetilde{G}$ as
above.
Hence, $G = P_1 \times P_2 \times \cdots \times P_c \times
\widetilde{G}$, and similarly,
$H = Q_1 \times Q_2 \times \cdots \times Q_d \times
\widetilde{H}$.  We say that $G$ is \textit{similar} to $H$
if, and only if, the following three conditions hold:
\begin{enumerate}
\item $\widetilde{G}$ is isomorphic to $\widetilde{H}$;
\item $c = d$;
\item $n_i = m_i$ for some reordering, where $|P_i|=p_i^{n_i}$
and $|Q_i| = q_i^{m_i}$.
\end{enumerate}
\end{definition}

From the comments before the definition, we see that if $G$ is
similar to $H$, then they will have the same number of subgroups.
Also, note that the equivalence class of $G$ is determined by
$\widetilde{G}$ and the multiset $[n_1, \ldots, n_c]$.
If we have a bound on the number of subgroups of $G$, then this
will bound $c$ and each $n_i$.  Therefore, the only remaining
hole in our theorem will be filled with the following lemma.

\begin{lemma}
Given $k > 1$ there are a finite number of
isomorphism classes of groups $G = \widetilde{G}$ (i.e. G
has no cyclic, central Sylow subgroups) having at most $k$
subgroups.
\end{lemma}
\end{quotation}

We wish to identify all similarity classes of groups with at most 12 subgroups.

To begin, we will find all groups with $G = \widetilde{G}$ and at most 12 subgroups.

\section{Abelian Groups}

If $G$ is abelian, then no sylow $p$-subgroup can be cyclic.  This means that if $G$ is
not a $p$-group, then it has at least $p+1$ subgroups of order $p$ and $q+1$ subgroups of
order $q$ for some distinct primes $p$ and $q$ and so at least $(p+1)(q+1) \geq 12$ non-trivial
subgroups.  So, $G$ must be a non-cyclic $p$-group.

If $G$ contains $Z_p \times Z_p \times Z_p$ then it has at least $2p^2 + 2p +4 \geq 16$ subgroups.
Consequently, $G = Z_{p^r} \times Z_{p^s}$ with $r \geq s \geq 1$.

Now, $Z_{p^2} \times Z_{p^2}$ has $p^2 + 3p + 5 \geq 15$ subgroups, so we must have $s = 1$.
Furthermore, $Z_{p^4} \times Z_p$ has $4p + 6 \geq 14$ subgroups, so $r \leq 3$.

Hence, the abelian groups with $G = \widetilde{G}$ and at most 12 subgroups are:

$r = 1$: ($p + 3$ subgroups) p = 2, 3, 5, 7

$r = 2$: ($2p + 4$ subgroups) p = 2, 3

$r = 3$: ($3p + 5$ subgroups) p = 2

\section{Non-Abelian Groups}

Now we consider non-abelian $G$.  Let $q$ be any prime dividing $|G|$.  If sylow $q$ is not normal, there
are at least $q+1$ sylow $q$-subgroups.  If sylow $q$ is normal, there is a $q$-complement.  If the complement
is normal, then the sylow $q$ cannot be cyclic (else $G \neq \widetilde{G}$) and so has at least $q+1$ maximal
subgroups. Finally, if the complement is not normal, then we have at least $q$ complements and 1 sylow $q$.
So, in every case, with the trival subgroup and $G$, the number of subgroups is at least $q+3$.  This implies
$q \leq 9$.

Can $|G|$ be divisible by all of 2, 3, 5, and 7?  If three of the sylows are not normal, counting conjugates of those sylows we have at least $(2 + 1) + (3 + 1) + (5 + 1) = 13$ subgroups, so at least 2 must be normal, say $P$ and $Q$.  Now their product $N$ has a complement $H$ with two primes and every subgroup of $H$ gives rise to at least four subgroups
of $G$ (itself, itself times $P$, itself times $Q$, and itself times $N$).  Consequently, $H$ has at most 3 subgroups and
so cannot have order divisible by two primes.

Suppose $|G|$ is divisible by 3 distinct primes $p, q, r$.
%[If sylow $p$ is normal, it must be central, and so non-cyclic.  So, at least $p + 1$ subgroups of
%order $p$.  If sylow $p$ not normal, must be at least $q$ sylow $p$-subgroups.]

If none of the sylows are normal, we have at least $(p + 1) + (q + 1) + (r + 1) \geq 13$ subgroups,
so at least one is normal. Call it an $r$-subgroup $N$. Complementing $N$, we have a two prime group $H$ in $G$
and every subgroup of $H$ gives rise to at least two subgroups of $G$ (itself and itself times $N$).  So, $H$ has at most 6 subgroups.

%% ???Must we have $H = \widetilde{H}$???  NO: $Z_6$ on $Z_7$.

\begin{lemma} A group $H$ divisible by 2 primes with at most 6 subgroups is either $Z_p \times Z_q,
Z_p \times Z_{q^2}$, or $S_3$.
\end{lemma}
\noindent
\begin{proof}
If $H$ is abelian, one sylow has 2 subgroups, the other 2 or 3 subgroups, so
$H = Z_p \times Z_q$ or $Z_p \times Z_{q^2}$.  If $H$ has a non-abelian sylow, the sylow has
at least 3 maximal subgroups, at least 1 prime order, 1 trivial, 1 whole subgroup, plus the other sylow $\Rightarrow \geq 7$
subgroups, which is too many.  So both sylows of $H$ are abelian and so at least one is not normal.
Assume WLOG $p < q$.  Then there are at least p+1 non-normal sylows, one other sylow, the whole
group, and trivial $\Rightarrow \geq p+4$ subgroups.  So we must have equality on our estimate which
means, $p = 2$,  sylow 2 is not normal, $p+1$ divides $|G|$, so $q = 3$, sylow 3 is normal, sylows
have no non-trivial subgroups $\Rightarrow H = S_3$.
\end{proof}

\vspace{0.2in}
If $H$ is abelian with cyclic sylows, then neither sylow can act trivially on $N$ or it would be
central cyclic in $G$.  Therefore $N$ cannot normalize the sylows or $H$ and so the number of
conjugates of each is at least $r$.  This gives at least $3r$ subgroups of $G$.  Furthermore we
have $G$, $N$, 1, and $N H_p$ and $N H_q$, so at least $3r + 5$ subgroups.  Consequently,
we must have $r = 2$.

Our list of subgroups includes only $N$ and 1 from subgroups of $N$.  So, if $N$ has at least 4
subgroups, then $G$ will have at least 13 subgroups.  Therefore, $N$ is cyclic of order 2 or 4.
But odd primes cannot act non-trivially on a cyclic group of order 2 or 4, so this situation is
impossible.

Last, we consider the case of $H = S_3$.  In this case, $r \geq 5$.
If $H$ acts trivially on $N$, then $N$ cannot be cyclic and so has at least 8 subgroups.  Since
$S_3$ has 6 subgroups, $G$ will have at least 48 subgroups.  So, $H$ acts non-trivially on
$N$ and a sylow 2, $H_2$, must also act non-trivally on $N$.  In particular, $N$ does not normalize $H$
or $H_2$ and, as above, this gives us at least $2r$ subgroups of $G$.  We also have $G$, $N$, and 1.
Since $r$ is at least 5, this implies $G$ has more than 12 subgroups.

%If neither sylow of $H$ is normal in $H$, we have at least $(p + 1) + (q + 1) \geq 7$ subgroups,
%so at least one sylow of $H$ is normal in $H$.  So $G$ has a sylow tower.  Also, the complement
%inside $H$ has at most 3 subgroups, so is cyclic $s$ or $s^2$ for some $s \in \{ p, q, r \}$.

%%%%%%%%%%%%%%%%%%%
Consequently, $|G|$ must be divisible by at most two primes.

\section{Groups with two primes in the order}

Suppose $|G| = p^a q^b$, $p < q$ primes and $G$ is non-abelian.

%If $G$ has order simply $pq$, then $G$ has $q + 3$ subgroups which is less than 12 since $q \leq 7$.
%This can happen with $(p,q)$ = (2,3), (2,5), (2,7), or (3,7).

The argument for two primes is more complicated than I would like.  I'll use indentation to organize
assumptions (like in computer code).
\begin{description}
\item[If neither sylow is cyclic,] then we have $(p+1)+(q+1)$ maximals in the sylows, the two sylows and the
whole group and trivial.  That is, the number of subgroups is at least $p + q + 6$.  This can only be less
than or equal to 12 when $p = 2, q = 3$.  Furthermore, the close estimate means each sylow must be normal
and prime squared order.  That forces $G$ to be abelian $\Rightarrow\Leftarrow$.

\item[If sylow $p$ is cyclic,]\quad \\

	\begin{description}
	\item[If sylow $p$ is normal,] then the sylow $q$ cannot act non-trivially and so the sylow $p$ is cyclic
and central $\Rightarrow\Leftarrow$.
	\item[else sylow $p$ is not normal,] then the number of sylow $p$ is $\geq q$.  \\
		\begin{description}
		\item[If sylow $q$ not cyclic,] we have at least $q+1$ maximals and so at least $2q+4$ subgroups.
Again the estimate forces sylow $q$ to be $Z_3 \times Z_3$ and $p = 2$, at least 10 subgroups.
Furthermore, we cannot have 9 sylow 2's, so there is a central subgroup of order 3.
Counting again, we have whole, trivial, one 9, four 3's, three sylow 2's, three sylow 2 times central 3 $\Rightarrow \geq 13$ subgroups. $\Rightarrow\Leftarrow$
		\item[else sylow $q$ is cyclic,] then $G$ must be supersolvable since both sylows are cyclic and it follows the sylow $q$
must be normal.  sylow $p$ must act non-trivially on $q$, so $(p,q) = (2,3), (2,5), (2,7)$, or $(3,7)$.\\
The only case where the action is not that of an element of order $p$ is when $p=2, q=5$ when $Z_4$ acts
faithfully on $Z_5$.  In this case, $G$ will have a quotient group of order 20 isomorphic to $Z_4$ acting
faithfully on $Z_5$.  Since this group has 14 subgroups, we do not need to consider such $G$.
		Therefore, we can assume the sylow $p$ acts as an element of order $p$.  The sylow $p$ fixes only the identity in sylow $q$ (only possible non-triv
action).  This means the sylow $p$ is self-normalizing and so the number of sylow $p$ subgroups is $q^b$.  Now the
subgroup of order $p^{a-1}$ centralizes the sylow $q$ and so is the intersection of the sylow $p$'s.  Hence the non-triv
$p$-subgroups are order $p, \ldots, p^{a-1}$, and $q^b$ subgroups order $p^a$.  There are $b$ non-triv $q$-subgroups
and $(a-1)\cdot b$ proper abelian subgroups divisible by $pq$.  Furthermore, for each non-triv, proper $q$-subgroup, order $q^c$, we can multiply that by the sylow $p$'s to get subgroups of order $p^a q^c$.  There will be $q^c$ different sylow $p$'s
giving the same subgroup of order $p^a q^c$ and so a total of $q^{b-c}$ such subgroups.
Consequently, the number of (non-abelian) proper subgroups
divisible by $p^a q$ will be $q + q^2 + \cdots + q^{b-1} = (q^b - q)/(q-1)$.
Including trivial and whole group, the number of subgroups is:
$$
(a-1 + q^b) + b + (a-1)b + \frac{q^b - q}{q-1} + 2 = q^b + \frac{q^b - q}{q-1} + a(b + 1) + 1.
$$
The corresponding values are shown in the following tables.

\qquad\begin{tabular}{|c c|c|c|}
\hline
\multicolumn{4}{|c|}{$q = 3$}\\
\hline
&\multicolumn{3}{c|}{$b$}\\
& &1&2\\
\cline{2-4}
&1&6&16\\
\cline{2-4}
$a$&2&8&19\\
\cline{2-4}
&3&10&22\\
\cline{2-4}
&4&12&25\\
\cline{2-4}
&5&14&28\\
\hline
\end{tabular}
\begin{tabular}{|c c|c|c|}
\hline
\multicolumn{4}{|c|}{$q = 5$}\\
\hline
&\multicolumn{3}{c|}{$b$}\\
& &1&2\\
\cline{2-4}
&1&8&34\\
\cline{2-4}
$a$&2&10&37\\
\cline{2-4}
&3&12&40\\
\cline{2-4}
&4&14&43\\
\hline
\end{tabular}
\begin{tabular}{|c c|c|c|}
\hline
\multicolumn{4}{|c|}{$q = 7$}\\
\hline
&\multicolumn{3}{c|}{$b$}\\
& &1&2\\
\cline{2-4}
&1&10&60\\
\cline{2-4}
$a$&2&12&63\\
\cline{2-4}
&3&14&66\\
\hline
\end{tabular}

Since the $q = 7$ table applies to both $p = 2$ and $p = 3$ we see there are 11 groups with at most 12
subgroups in this case.\\

		\end{description}

	\end{description}
	\item[else sylow $p$ is not cyclic,] and so sylow $q$ must be cyclic and not central.

	\begin{description}
	\item[If neither sylow is normal,] then number of sylow $p$ is at least $q$ and the
number of sylow $q$ is at least $q+1$.  Including trivial and whole, we have at least
$2q + 3$ subgroups, which means $q < 5$ and so we have $q = 3, p = 2$.
We have at least three sylow 2's and at least four sylow 3's, triv, whole, and one
sylow 2 will contain three subgroups of order 2.  However, that's already 12 subgroups
and we have more order 2 subgroups in the other sylow 2's. $\Rightarrow\Leftarrow$
	\item[If sylow $p$ is normal,] then it cannot be cyclic, else it would be central.  Furthermore,
$Z_p \times Z_p \times Z_p$ has too many subgroups, so the sylow $p$ must be a two-generator
group.  So, the sylow $p$ has $p+1$ maximal subgroups, there are at least $q+1$ sylow $q$'s
and with sylow $p$ itself, triv, and whole, we have at least $p + q + 5$ subgroups.  If $p \geq 3$,
this is too many, so we must have $p = 2$.  If $q \geq 5$, then our count gives at least 12 subgroups.
The only way to avoid going over 12 would be to have the sylow $p$ be $Z_2 \times Z_2$, sylow
$q$ be $Z_5$.  But $Z_5$ has no non-trivial action on $Z_2 \times Z_2$, and so we are left
only with the case $p = 2, q = 3$, and at least 10 subgroups.  Now, the alternating group
$A_4$ satisfies these conditions and has 10 subgroups.  If sylow 3 has order larger than 3, then
there will be a central 3-subgroup whose product with the various 2-subgroups will give more
than 12 subgroups.  So, the sylow 3 must be order 3.  If sylow 2 has order larger than 4, then
the frattini subgroup of the sylow 2 and the frattini subgroup times the various sylow 3-subgroups give more
than 12 subgroups.  So, $A_4$ is the only group in this case.
	\item[If sylow $q$ is normal,] we are already assuming the sylow $q$ is cyclic.
Furthermore, the sylow $p$ is not normal and not cyclic.  So we have at least $q$ sylow $p$'s and
the sylow $p$'s have at least $p+1$ maximal subgroups.
		\begin{description}
		\item[If two sylow $p$'s have a common maximal subgroup, $M$,] then the normalizer of $M$ will contain
at least two sylow $p$-subgroups of $G$ and so must be divisible by $q$.  Thus the intersection of the normalizer with
the sylow $q$-subgroup is a non-trivial central $q$ subgroup.
Therefore, our subgroups include $q$ sylow $p$'s, $q$ sylow $p$'s times the central $q$, at least $p+1$ maximal
subgroups of sylow $p$, each of those $p+1$ times the central $q$, 
the sylow $q$, the central $q$, triv, and whole for at least $2q +2p + 6 \geq 16$ subgroups. $\Rightarrow\Leftarrow$
		\item[else no shared maximal subgroups,] which means we have at least $q(p+1)$ maximal subgroups of sylow $p$'s
altogether.  Therefore we have at least  these $q(p+1)$ subgroups, $q$ sylow $p$-subgroups, one sylow $q$, triv, and whole.
That is, $qp +2q + 3 \geq 15$ subgroups. $\Rightarrow\Leftarrow$

		\end{description}

	\end{description}
\end{description}

Hence, the non-abelian groups with $G = \widetilde{G}$, $|G|$ divisible by at least two primes, and at most 12 subgroups are:

$Z_p \ltimes Z_q$: ($q+3$ subgroups) $(p, q)$ = (2, 3), (2, 5), (2, 7), or (3, 7)

$Z_{p^2} \ltimes Z_q$ with action of order $p$: ($q+5$ subgroups) $(p, q)$ = (2, 3), (2, 5), (2, 7), or (3, 7)

$Z_8 \ltimes Z_q$ with action of order 2: ($q+7$ subgroups) $q$ = 3 or 5

$Z_{16} \ltimes Z_3$ with action of order 2: (12 subgroups)

$A_4$: (10 subgroups)

%%%%%%%%%%%%%%%%%%%
\section{$p$-Groups}

Finally, we have the case of $G$ a non-abelian $p$-group, $|G| = p^n$.

%Known cases:
%
%$D_8$:(10), $Q_8$:(6), gen'l quaternion order 16:(11),
%$Z_2 \ltimes Z_8$:(11), extraspecial 27, exponent 9:(10)

The group $Z_p \times Z_p \times Z_p$ has  at least 16 subgroups, so $G$ must be a two-generator group.
Consequently, we have $G$, $p+1$ maximals, the Frattini subgroup, and triv $\Rightarrow \geq p+4$ subgroups.

When $p$ is odd we will have at least $p+1$ order $p$, one of which might be the Frattini, so we'd have
at least $2p+4 \Rightarrow p$ = 2 or 3.

Assume $|G'| = p$ and $Z(G)$ cyclic.  Then $G' \subset Z(G)$ and for $x,y \in G$, $[x^p,y] = [x,y]^p = 1$
so $x^p$ is central and $\Phi(G) \subset Z(G)$.  Now $|G: \Phi(G)| = p^2$ and so we must have $\Phi(G) = Z(G)$
is cyclic order $p^{n-2}$.

If $G$ has a single subgroup of order $p$, then $p = 2$ and $G$ is generalized quaternion.  Now, each generalized
quaternion 2-group contains the next smaller as a subgroup.  So, we can just check: $Q_8$ has 6 subgroups,
$Q_{16}$ has 11 subgroups, and $Q_{32}$ has 20.  That is, $G$ must be $Q_8$ or $Q_{16}$.

Otherwise, $G$ has more than one subgroup of order $p$ and so we can choose $S \subset G$ with $S$ order p and
not in the center.  Then $M = S Z(G)$ is an abelian subgroup isomorphic to $Z_p \times Z_{p^{n-2}}$.  From the
abelian case above, we see that the number of subgroups of $M$ is given in the following table.

\qquad\begin{tabular}{|c c|c|c|c|} 
\hline
&\multicolumn{4}{c|}{$n$}\\
& &3&4&5\\
\cline{2-5}
$p$&2&5&8&11\\
\cline{2-5}
&3&6&10&14\\
\hline
\end{tabular}

In addition to the subgroups in $M$, $G$ also has $p$ other maximal subgroups, and $G$ itself.  So, increasing
each table entry by $p+1$ we see the only possibilities for at most 12 subgroups are $|G|$ = 8, 16, or 27.

By brute force check we find for order 8 the dihedral group $D_8$ with 10 subgroups, for order 16 a group with presentation $\langle a, b|a^2, b^8, b^a = b^5 \rangle$
having 11 subgroups, and
for order 27 the extraspecial group of exponent 9 with 10 subgroups.

Thus we have found five groups with $|G'| = p$ and cyclic center.  Now any non-abelian finite $p$-group will
have such a group as a homomorphic image.

%Let's consider when $G$ might have $Q_8$ as a homomorphic image.
%
%Let $Q = G/G'$.  Then $Q$ is a two generator abelian 2-group.  From above we see that $|Q| \leq 16$.

Note that no generalized quaternion group can have a smaller generalized quaternion group as a homomorphic image.
One easy way to see this is to note that $Q_n/Z(Q_n) = D_{n/2}$ and any non-abelian image of a dihedral group will have more
than one involution.  We will use this fact several times below.

Suppose $|G| = 3^5$ and choose $K$ maximal such that $\overline{G} = G/K$ is non-abelian.  It follows that
$\overline{G}$ will have $|\overline{G}'| = 3$ and $Z(\overline{G})$ cyclic.  If $|K| \leq 3$, it follows from
above that $\overline{G}$, and so $G$, will have more than 12 subgroups.  If $|K| = 3^2$, then $|\overline{G}| = 3^3$
and so must be the extraspecial group of exponent $3^2$.  Thus, $\overline{G}$ has 10 subgroups.  The only way
we could have at most 12 subgroups in $G$ is if $K$ is cyclic and every subgroup either contains $K$ or is order 3
or 1 in $K$.  However, this implies that $G$ has only one subgroup of order 3, which is impossible.  Consequently, no
non-abelian 3-group of order at least $3^5$ can have at most 12 subgroups.  We know from above that there is
only one such group of order $3^3$ and a computer check shows that there are no examples of order $3^4$ (these
groups have at least 14 subgroups).

Now suppose $|G| = 2^6$ and choose $K$ maximal such that $\overline{G} = G/K$ is non-abelian.  As above,
$\overline{G}$ will have $|\overline{G}'| = 2$ and $Z(\overline{G})$ cyclic.  If $|K| \leq 2$, it follows from
above that $\overline{G}$, and so $G$, will have more than 12 subgroups.
%If $|K| = 2$, then $\overline{G}$ must be $Q_{16}$ with 11 subgroups.  If every subgroup of $G$ either contains $K$
%or is contained in $K$, then $G$ would have 12 subgroups.  However, this would force $G$ to have a single
%subgroup of order 2, making it $Q_{32}$ which it is not (and would still be too many subgroups even if it were).
  If $|K| = 2^2$, then $|\overline{G}| = 2^4$
and so we see from above that $\overline{G}$ is one of two groups each of which have 11 subgroups.  Since $K$
has at least 2 proper subgroups, $G$ must have 13 or more subgroups.  Finally, consider $|K| = 2^3$.  Then $\overline{G}$
must be $D_8$ or $Q_8$ with 10 or 6 subgroups respectively.  Since $K$ has at least 3 proper subgroups, $\overline{G}$ cannot be $D_8$ and so must be $Q_8$ with 6 subgroups.  Even if all of the subgroups of $G$ either contained $K$ or were
contained in $K$, then $K$ would have to have at most 7 subgroups (6 proper).  However, the only groups of order $2^3$
with 7 or fewer subgroups are $Z_8$ and $Q_8$.  Since each of these have a single subgroup of order 2, our
assumption would force $G$ to be generalized quaternion, which it clearly is not.  Thus we see that no non-abelian 2-group
of order at least $2^6$ can have number of subgroups less than or equal to 12.  A computer check shows there are
no examples of order $2^5$ (these groups have at least 14 subgroups) and only the two groups mentioned above for
order $2^4$.

Hence, the non-abelian $p$-groups with $G = \widetilde{G}$ and at most 12 subgroups are:

$D_8$: (10 subgroups)

$Q_8$: (6 subgroups)

$Q_{16}$: (11 subgroups)

$Z_2 \ltimes Z_8$: (11 subgroups)

$E_{27}$ extraspecial order 27, exponent 9: (10 subgroups)

\section{Conclusion}

Collecting all of our results, we have:

\qquad\begin{tabular}{c|l}
$n$&Groups with $G = \widetilde{G}$ and $n$ subgroups\\
\hline
1&Trivial group\\
2&\\
3&\\
4&\\
5&$Z_2 \times Z_2$\\
6&$Z_3 \times Z_3$, $S_3$, $Q_8$\\
7&\\
8&$Z_2 \times Z_4$, $Z_5 \times Z_5$, $D_{10}$, $Z_4 \ltimes Z_3$\\
9&\\
10&$Z_3 \times Z_9$, $Z_7 \times Z_7$, $D_{14}$, $A_4$, $Z_3 \ltimes Z_7$, $Z_4 \ltimes Z_5$, $Z_8 \ltimes Z_3$, $D_8$, $E_{27}$\\
11&$Z_2 \times Z_8$, $Q_{16}$, $Z_2 \ltimes Z_8$\\
12&$Z_4 \ltimes Z_7$, $Z_9 \ltimes Z_7$, $Z_8 \ltimes Z_5$, $Z_{16} \ltimes Z_3$\\
\end{tabular}

In particular the sequence of the number of groups with $G = \widetilde{G}$ and $n$ subgroups would be:
$$
1, 0, 0, 0, 1, 3, 0, 4, 0, 9, 3, 4, \ldots
$$

Forming the direct product with a coprime, cyclic group of order $p^k$ will multiply the number of subgroups by $k+1$.
Thus we find groups with $n$ subgroups corresponding to various factorizations of $n$.

\begin{tabular}{c|l}
$n$&Similarity class representatives with $n$ subgroups\\
\hline
1&Trivial group\\
2&$Z_2$\\
3&$Z_{2^2}$\\
4&$Z_{2^3}$, $Z_2 \times Z_3$\\
5&$Z_{2^4}$, $Z_2 \times Z_2$\\
6&$Z_{2^5}$, $Z_2 \times Z_{3^2}$, $Z_3 \times Z_3$, $S_3$, $Q_8$\\
7&$Z_{2^6}$\\
8&$Z_{2^7}$, $Z_2 \times Z_{3^3}$, $Z_2 \times Z_3 \times Z_5$, $Z_2 \times Z_4$, $Z_5 \times Z_5$, $D_{10}$, $Z_4 \ltimes Z_3$\\
9&$Z_{2^8}$, $Z_{2^2} \times Z_{3^2}$\\
10&$Z_{2^9}$, $Z_2 \times Z_{3^4}$, $Z_2 \times Z_2 \times Z_3$, $Z_3 \times Z_9$, $Z_7 \times Z_7$, $D_{14}$, $A_4$, $Z_3 \ltimes Z_7$, $Z_4 \ltimes Z_5$,\\
&\qquad $Z_8 \ltimes Z_3$, $D_8$, $E_{27}$\\
11&$Z_{2^{10}}$, $Z_2 \times Z_8$, $Q_{16}$, $Z_2 \ltimes Z_8$\\
12&$Z_{2^{11}}$, $Z_2 \times Z_{3^5}$, $Z_3 \times Z_3 \times Z_2$, $S_3 \times Z_5$, $Q_8 \times Z_3$, 
$Z_{2^2} \times Z_{3^3}$, $Z_2 \times Z_3 \times Z_{5^2}$,\\
&\qquad $Z_4 \ltimes Z_7$, $Z_9 \ltimes Z_7$, $Z_8 \ltimes Z_5$, $Z_{16} \ltimes Z_3$\\
\end{tabular}

In the following version of the previous table we represent classes using $p, q, r$ to represent primes which do not
occur anywhere else in the group order.  This makes it a bit easier to recognize the infinite classes.

\begin{tabular}{c|l}
$n$&Similarity classes with $n$ subgroups\\
\hline
1&Trivial group\\
2&$Z_p$\\
3&$Z_{p^2}$\\
4&$Z_{p^3}$, $Z_p \times Z_q$\\
5&$Z_{p^4}$, $Z_2 \times Z_2$\\
6&$Z_{p^5}$, $Z_p \times Z_{q^2}$, $Z_3 \times Z_3$, $S_3$, $Q_8$\\
7&$Z_{p^6}$\\
8&$Z_{p^7}$, $Z_p \times Z_{q^3}$, $Z_p \times Z_q \times Z_r$, $Z_2 \times Z_4$, $Z_5 \times Z_5$, $D_{10}$, $Z_4 \ltimes Z_3$\\
9&$Z_{p^8}$, $Z_{p^2} \times Z_{q^2}$\\
10&$Z_{p^9}$, $Z_p \times Z_{q^4}$, $Z_2 \times Z_2 \times Z_p$, $Z_3 \times Z_9$, $Z_7 \times Z_7$, $D_{14}$, $A_4$, $Z_3 \ltimes Z_7$, $Z_4 \ltimes Z_5$,\\
&\qquad $Z_8 \ltimes Z_3$, $D_8$, $E_{27}$\\
11&$Z_{p^{10}}$, $Z_2 \times Z_8$, $Q_{16}$, $Z_2 \ltimes Z_8$\\
12&$Z_{p^{11}}$, $Z_p \times Z_{q^5}$, $Z_3 \times Z_3 \times Z_p$, $S_3 \times Z_p$, $Q_8 \times Z_p$, 
$Z_{p^2} \times Z_{q^3}$, $Z_p \times Z_q \times Z_{r^2}$,\\
&\qquad $Z_4 \ltimes Z_7$, $Z_9 \ltimes Z_7$, $Z_8 \ltimes Z_5$, $Z_{16} \ltimes Z_3$\\
\end{tabular}

In conclusion, the sequence of number of similarity classes with a given number of subgroups begins:
$$
1, 1, 1, 2, 2, 5, 1, 7, 2, 12, 4, 11, \ldots
$$

Note: In \cite{Miller}, Miller lists the groups with specified number of subgroups where the number of subgroups runs from 1 through 9.  His lists agree with ours.

\end{document}